\documentclass[a4paper,11pt]{article}

\usepackage{fourier}
\usepackage{color}
\usepackage{graphicx}
\usepackage{url}
\usepackage[affil-it]{authblk}
\usepackage{amsmath}
\usepackage{wrapfig}
\usepackage{xspace}

\usepackage[T1]{fontenc}
\usepackage{times}

\usepackage{subcaption}
\usepackage{multirow}

\begin{document}

\title{A Restricted-Domain Dual Formulation for Two-Phase Image Segmentation}

\author{Jack Spencer}
\affil{Department of Mathematics, \\ University of Liverpool, UK.}
\date{}
\maketitle
\thispagestyle{empty}

\begin{abstract}
In two-phase image segmentation, convex relaxation has allowed global minimisers to be computed for a variety of data fitting terms. Many efficient approaches exist to compute a solution quickly. However, we consider whether the nature of the data fitting in this formulation allows for reasonable assumptions to be made about the solution that can improve the computational performance further. In particular, we employ a well known dual formulation of this problem and solve the corresponding equations in a restricted domain. We present experimental results that explore the dependence of the solution on this restriction and quantify imrovements in the computational performance. This approach can be extended to analogous methods simply and could provide an efficient alternative for problems of this type.
\end{abstract}

\section{Introduction}

Image segmentation is the meaningful partitioning of an image based on certain characteristics. In two-phase segmentation this consists of determining the foreground and background of a domain $\Omega\in\mathbb{R}^{2}$, i.e. find a closed boundary separating subregions $\Omega_{1}$ and $\Omega_{2}$. This is distinct from multiphase approaches, where more than two separate regions are determined. Our work concerns the continuous setting, which we will briefly discuss in the next section. Equivalent problems in the discrete setting have been well studied with details found in \cite{BoyKol:04}. A comprehensive background behind the following functional can also be found in \cite{CPintro}. Briefly, the aim is to determine an indicator function, $u(x)$, that labels the foreground and background by minimising the following energy:
\begin{equation} \label{eqn:F}
\min_{u\in\lbrace0,1\rbrace}\left\lbrace\int_{\Omega}|\nabla u(x)|\ dx+\lambda\int_{\Omega}f(x)u(x) dx\right\rbrace.
\end{equation}
The function $f(x)$ is typically referred to as the fitting term, determining how the segmentation solution corresponds to the data. It is balanced by a regularisation term, in this case the total variation (TV) semi-norm which penalises the length of the segmentation boundary. For large $\lambda$ the following will hold precisely:
\[
\begin{array}
[c]{l}%
f(x)<0,\ \ \ \text{foreground}, \\
f(x)>0,\ \ \ \text{background}.
\end{array}
\]
When $\lambda$ is varied the solution for $u$ will fit the data with more regularity, i.e. some areas where $f<0$ will be background and some areas where $f>0$ will be foreground. However, this is most likely where $f$ is close to 0. With that in mind our work considers what improvements can be made by concentrating on regions in the domain where $f(x)\approx0$. An exception to this concerns cases where $f$ has strong noise, which we will return to later. A number of choices for $f$ exist depending on the application, including piecewise-constant segmentation based on the work of \cite{ACWE}]:
\begin{equation} \label{eqn:CV}
f(x)=(z-c_{1})^2-(z-c_{2})^2,
\end{equation}
where $c_{1}$ and $c_{2}$ are intensity constants indicating average foregound and background intensitites of the image $z(x)$, respectively. We also consider a selection based fitting term such as \cite{CDSS}
\begin{equation} \label{eqn:DS}
f(x)=(z-c_{1})^2-(z-c_{2})^2+\gamma P(x),
\end{equation}
where $\gamma P(x)$ is a distance selective term based on user input. In Section \ref{ER} we present results using equations \eqref{eqn:CV} and \eqref{eqn:DS}. This approach is not limited to the fitting terms mentioned above, and can be extended to any segmentation problem  in this framework. Alternatives for future consideration include bias field segmentation \cite{DC} and interactive convex active contours \cite{Nguyen:12}.

A restricted-domain approach is analogous to banded segmentation methods such as \cite{Rommelse:03} and \cite{Zhang:14}, among many others. This work differs in the sense that the restriction is based on the values of the fitting term, rather than the location of the boundary at an iteration, which is potentially simpler computationally. We compute an approximation of the global minimiser of the energy \eqref{eqn:F}, with the accuracy determined by the level of domain restriction. This is based on the following initialisation of the indicator function: $u^{(0)}=H(-f)$, where $H(\cdot)$ is the Heaviside function. 

In the following we will briefly introduce existing methods for finding the global minimisers of two-phase segmentation problems, introducing the dual formulation of \cite{Bresson:07} based on \cite{Chambolle:04}. We then discuss the proposed approach where a restricted domain based on the fitting term is considered, before detailing the method and how it relates to \cite{Bresson:07}. Finally, we present some results for three examples for various restrictions on the domain, quantifying the accuracy and computational performance in comparison to the original method. We then offer some concluding remarks.

\section{Convex Relaxation for Two-Phase Segmentation}

We now introduce the details of the approach we consider in this work. Again, a comprehensive background of this work can be found in \cite{CPintro} and many others. Essentially, convex relaxation in this case involves relaxing the binary constraint in the original functional \eqref{eqn:F}, i.e. $u\in[0,1]$. The seminal work here is \cite{Chan:06} who found global minimisers of the two-phase piecewise-constant Mumford-Shah model \cite{MumfordShah} (assuming fixed intensity constants, $c_{1}, c_{2}$). Therefore, the problem considered here is:
\begin{equation} \label{eqn:CR}
\min_{u\in[0,1]}\left\lbrace\int_{\Omega}|\nabla u(x)|\ dx+\lambda\int_{\Omega}f(x)u(x) dx\right\rbrace.
\end{equation}
In \cite{Chan:06} they introduce a penalty function to enforce the constraint on $u$ and solve using time marching. It is also possible to use additive operator splitting \cite{CDSS}, split Bregman \\ \cite{Goldstein:10}, and Chambolle-Pock \cite{ChambollePock} among many others. However, initially we intend to implement our restricted-domain approach on the dual formulation used in \\ \cite{Bresson:07}. We breifly detail this approach next.

\subsection{Dual Formulation} \label{CDF}

The dual formulation of this problem was first introduced by \cite{Bresson:07}, based on the work of \cite{Chambolle:04}, \cite{Aujol:06} and the references therein. The idea is to introduce a new variable, $v(x)$, and minimise the following functional alternately:
\begin{equation} \label{eqn:DF}
\min_{u,v}\left\lbrace\int_{\Omega}|\nabla u(x)|\ dx+\frac{1}{2\theta}\int_{\Omega}\left(u(x)-v(x)\right)^{2}\ dx+\int_{\Omega}\lambda f(x)v(x)+\alpha\psi(v)\ dx\right\rbrace,
\end{equation}
where $\psi(v)=\max\lbrace0,2|v-1/2|-1\rbrace$. By splitting the variables in this way, the minimisation of $u$ concentrates on the TV term, and the minimisation of $v$ satisifes the fitting and constraint requirements. In \cite{Bresson:07} the regularisation term is weighted, however here we concentrate on the original problem \eqref{eqn:CR}. The parameter $\alpha$ ensures the constraints on the indicator function $u(x)$ in \eqref{eqn:DF} are met, and can be set automatically \cite{Chan:06}. The minimistion of $u$ and $v$ can be achieved iteratively by the following steps. With fixed $v$, the solution of $u$ is given by
\begin{equation*}
u(x)=v(x)-\theta\nabla\cdot\rho(x),
\end{equation*}
where $\rho=(\rho^{1},\rho^{2})$ is the solution of
\begin{equation*}
\nabla(\theta\nabla\cdot\rho-v)-|\nabla(\theta\nabla\cdot\rho-v)|\rho=0,
\end{equation*}
which can be solved by a fixed point method. With fixed $u$, the solution for $v$ is given as
\begin{equation*}
v(x)=\min\lbrace\max\lbrace u(x)-\theta\lambda r(x),0\rbrace,1\rbrace.
\end{equation*}
This is repeated until convergence. Further details can be found in \cite{Bresson:07}, including the definition of the discrete gradient and divergence operators from \cite{Chambolle:04}. 

\section{Proposed Approach: Segmentation in a Restricted Domain}

\begin{figure}[!ht]
\begin{center} 
    \begin{subfigure}{.29\linewidth}
\includegraphics[width=\linewidth]{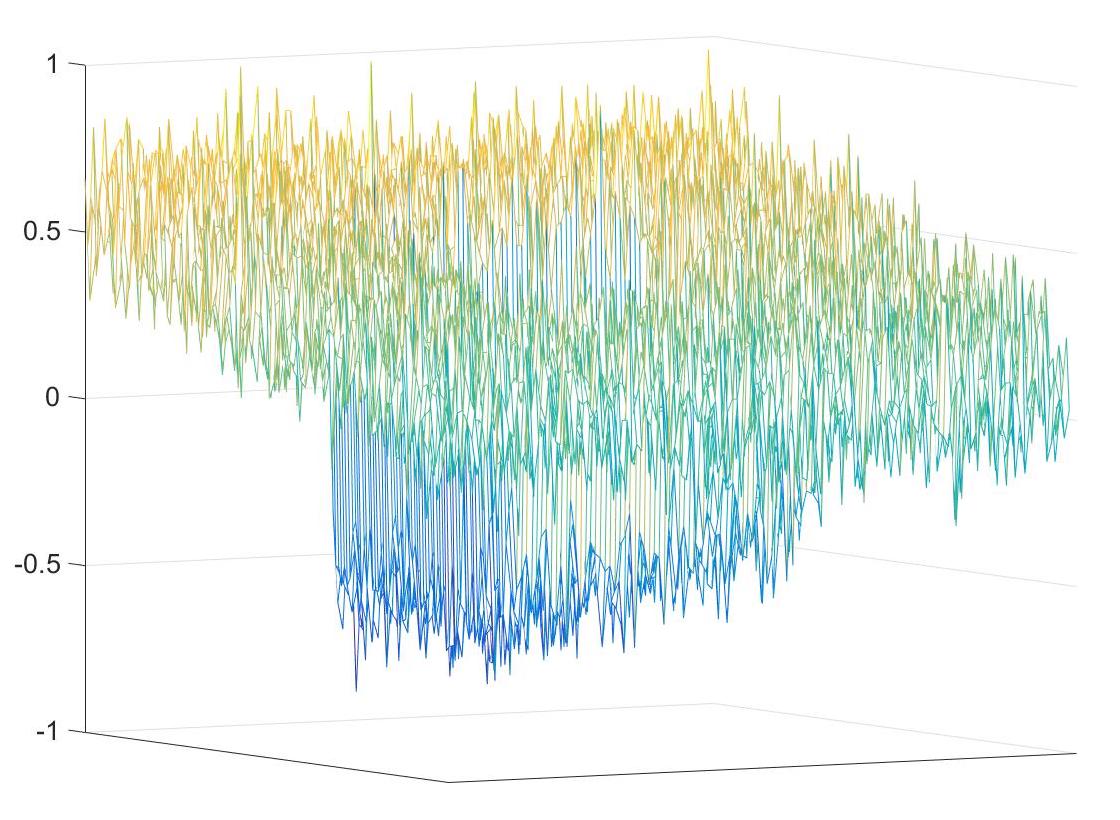}
    \end{subfigure}
    \begin{subfigure}{.29\linewidth}
\includegraphics[width=\linewidth]{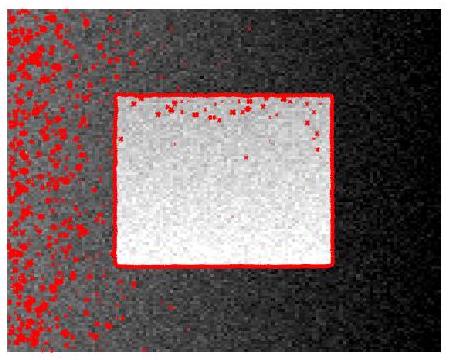}
    \end{subfigure}

   \begin{subfigure}{.29\linewidth}
\includegraphics[width=\linewidth]{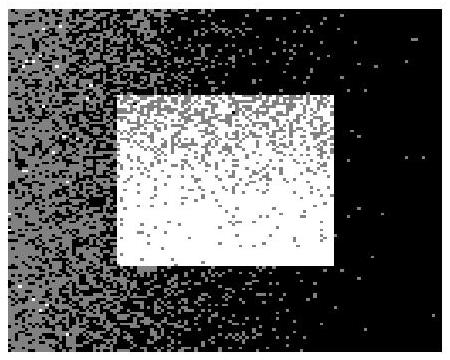}
    \end{subfigure}     
   \begin{subfigure}{.29\linewidth}
\includegraphics[width=\linewidth]{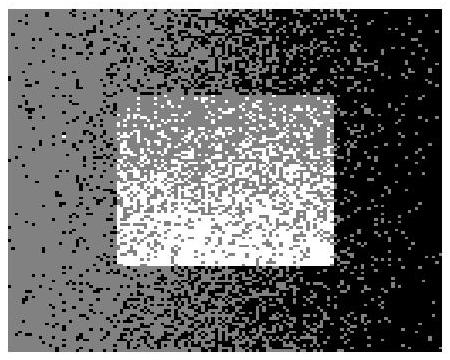}
   \end{subfigure}
    \begin{subfigure}{.29\linewidth}
\includegraphics[width=\linewidth]{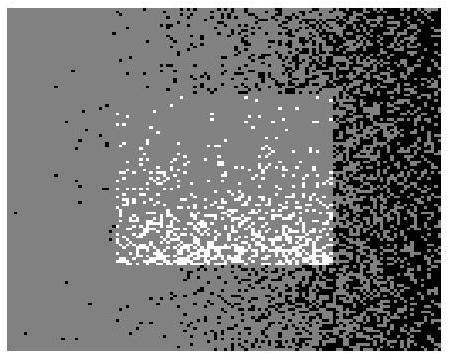}
    \end{subfigure} 
\end{center}
\vspace{-20pt}
\caption{Restricted Domain: The top row shows the fitting term (left) and the image (right) with the zero level-set of $f(x)$ given in red. The bottom row indicates the corresponding regions $Fg$ (white), $Bg$ (black), and $RD$ (grey) for values of $q=0.25,0.5,0.75$ from left to right, respectively.} \label{fig:myfig}
\vspace{-10pt}
\end{figure}

We now introduce our approach for reducing the computation time for this type of problem. We begin by assuming the solution in certain parts of the discretised domain based on the values of $f(x)$ for a given problem. The indicator function is then fixed at 1 or 0 at these points for the foreground ($Fg$) and background ($Bg$) respectively. We solve the equation in the remaining region, which we call the restricted domain ($RD$). Let us define $q\in[0,1]$ such that the following thresholding of the fitting function holds. We define a value $\hat{q}\in\mathbb{R}$ such that the percentage of nodes in the restricted domain is $q\ (\times100)$:
\[
\left\{
\begin{array}{cl}
Fg&=\ \left\lbrace x:x\in\Omega,\ f(x)\leq-\hat{q}\right\rbrace \\
Bg&=\ \left\lbrace x:x\in\Omega,\ f(x)\geq\hat{q}\right\rbrace \\
RD&=\ \left\lbrace x:x\in\Omega\setminus Fg\setminus Bg\right\rbrace.
\end{array}
\right.
\]
The value of $\hat{q}$ is initially 0 and is increased until the selected value of $q$ is satisfied. In other words, for $q=0$, $RD=\emptyset$ and the solution is a zero-thresholding of the fitting function, or equivalent to selecting a large $\lambda$ in the original problem \eqref{eqn:CR}. For $q=1$, $Fg=\emptyset$ and $Bg=\emptyset$ and we consider the problem in a conventional manner with no restriction on the domain. For $q\in(0,1)$ we consider a restricted domain of varying degrees as illustrated in Figure \ref{fig:myfig}, where the corresponding region of interest is given in grey. This means that we need to minimise the energy in a restricted domain which, when combined with the dual formulation of \cite{Bresson:07} and \cite{Chambolle:04}, means that we are solving the following equation for $\rho$ at certain points:
\begin{equation} \label{eqn:rho}
\nabla(\theta\nabla\cdot\rho-v)-|\nabla(\theta\nabla\cdot\rho-v|\rho=0.
\end{equation}

This involves making certain assumptions about the solution for $\rho$. In \cite{Bresson:07}, it is initialised as $\rho^{(0)}=0$ and the solution of \eqref{eqn:rho} is clearly dependent on $u$ and $v$. However, for the initialisation $u^{(0)}=-H(f)$ the solution of $\rho$ has a predictable form. In particular, $\rho\approx0$ where $|f|$ is largest and $\rho\in[-1,1]$ where $|f|$ is closer to 0. If $q$ is selected sensibly (we will return to this later), when $u^{(0)}=-H(f)$ it is reasonable to assume that the solution of \eqref{eqn:rho} for $x\in Fg\cap Bg$ is $\rho^*(x)=0$. Clearly, the larger $q$ is the less reliable this assumption becomes and the corresponding solution for $\rho$ will be less accurate. However, part of this work concerns what consists of a sensible selection for $q$ and whether it is possible to make reasonable assumptions about $f$ that can improve the efficiency of minimising the original formulation \eqref{eqn:CR}. 

We now ellaborate on the details behind minimising \eqref{eqn:CR} with a dual formulation in a restricted domain. We first consider the following minimisation problem:
\begin{equation*}
\min_{u}\left\lbrace\int_{\Omega}|\nabla u(x)|\ dx+\frac{1}{2\theta}\int_{\Omega}\left(u(x)-v(x)\right)^{2}\ dx\right\rbrace.
\end{equation*}
The solution, based on our approach, is given by 
\begin{equation}
u(x)=\left\{
\begin{array}{ll}
1, &\text{for}\ x\in Fg \\
0, &\text{for}\ x\in Bg \\
v-\theta\nabla\cdot\rho, &\text{for}\ x\in RD,
\end{array}
\right.
\end{equation}
where $\rho=(\rho_{1},\rho_{2})$ satisfies 
\begin{equation}
\rho=\left\{
\begin{array}{ll}
0, &\text{for}\ x\in Fg\cap Bg \\
\nabla(\theta\nabla\cdot\rho-v)-|\nabla(\theta\nabla\cdot\rho-v)|\rho=0, &\text{for}\ x\in RD.
\end{array}
\right.
\end{equation}
For $x\in RD$ the following fixed point method, with time step $\tau$, will solve the equation for $\rho$:
\begin{equation*}
\rho^{n+1}=\frac{\rho^{n}+\tau\nabla(\nabla\cdot \rho^{n}-v/\theta)}{1+\tau|\nabla(\nabla\cdot\rho^{n}-v/\theta)|}
\end{equation*}
As before, the following minimisation problem is then solved with $u$ fixed:
\begin{equation*}
\min_{v}\left\lbrace\frac{1}{2\theta}\int_{\Omega}\left(u(x)-v(x)\right)^{2}\ dx+\int_{\Omega}\lambda f(x)v(x)+\alpha\psi(v)\ dx\right\rbrace.
\end{equation*}
We combine our assumptions about $u$ and $\rho$ with the work of \cite{Bresson:07} to give the corresponding solution as 
\begin{equation}
v(x)=\left\{
\begin{array}{ll}
1, &\text{for}\ x\in Fg \\
0, &\text{for}\ x\in Bg \\
\min\lbrace\max\lbrace u(x)-\theta\lambda f(x),0\rbrace,1\rbrace, &\text{for}\ x\in RD.
\end{array}
\right.
\end{equation}
As with \cite{Bresson:07}, as discussed in the previous section, $u$ and $v$ are minimised alternately until convergence. The main advantage of this approach concerns finding the solution of $\rho$ at each iteration with the fixed point method detailed above. Based on the choice of $q$ it is possible that significant advantages exist in terms of computation time with minimal compromise on the quality of the solution. We will discuss some exceptions to this, as well as future considerations in the following sections.

\section{Experimental Results} \label{ER}

\begin{figure}[!ht]
\begin{center} 
    \begin{subfigure}{.29\linewidth}
\includegraphics[width=\linewidth]{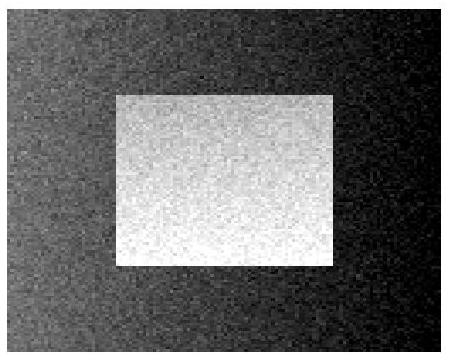}
   \end{subfigure}
    \begin{subfigure}{.29\linewidth}
\includegraphics[width=\linewidth]{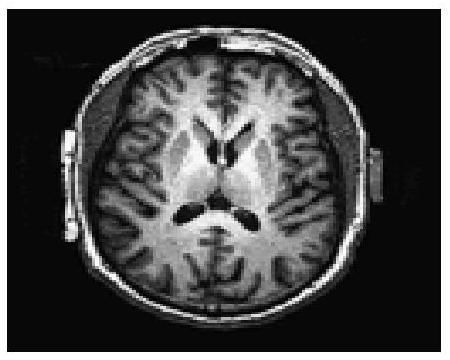}
    \end{subfigure}
    \begin{subfigure}{.29\linewidth}
\includegraphics[width=\linewidth]{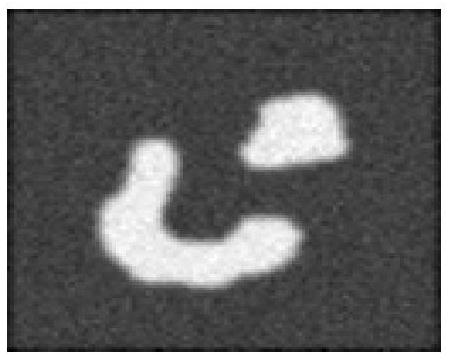}
    \end{subfigure}

    \begin{subfigure}{.29\linewidth}
\includegraphics[width=\linewidth]{h1.jpg}
   \end{subfigure}
    \begin{subfigure}{.29\linewidth}
\includegraphics[width=\linewidth]{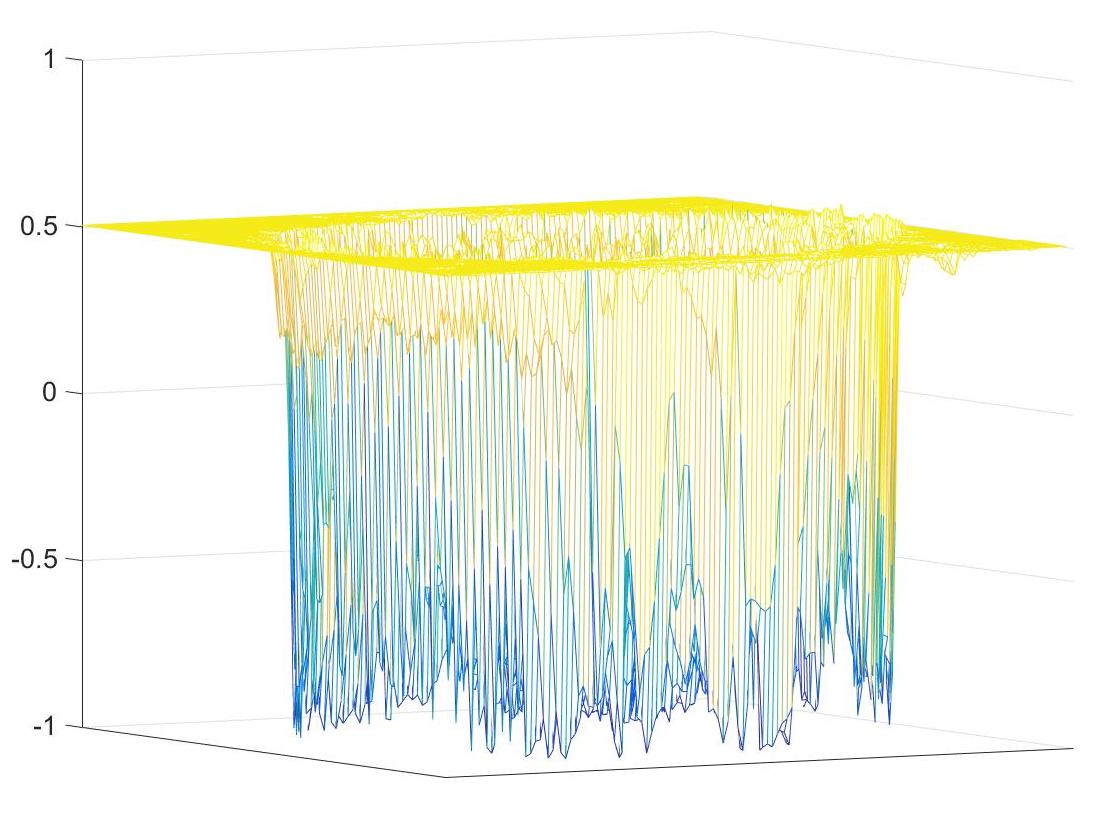}
    \end{subfigure}
    \begin{subfigure}{.29\linewidth}
\includegraphics[width=\linewidth]{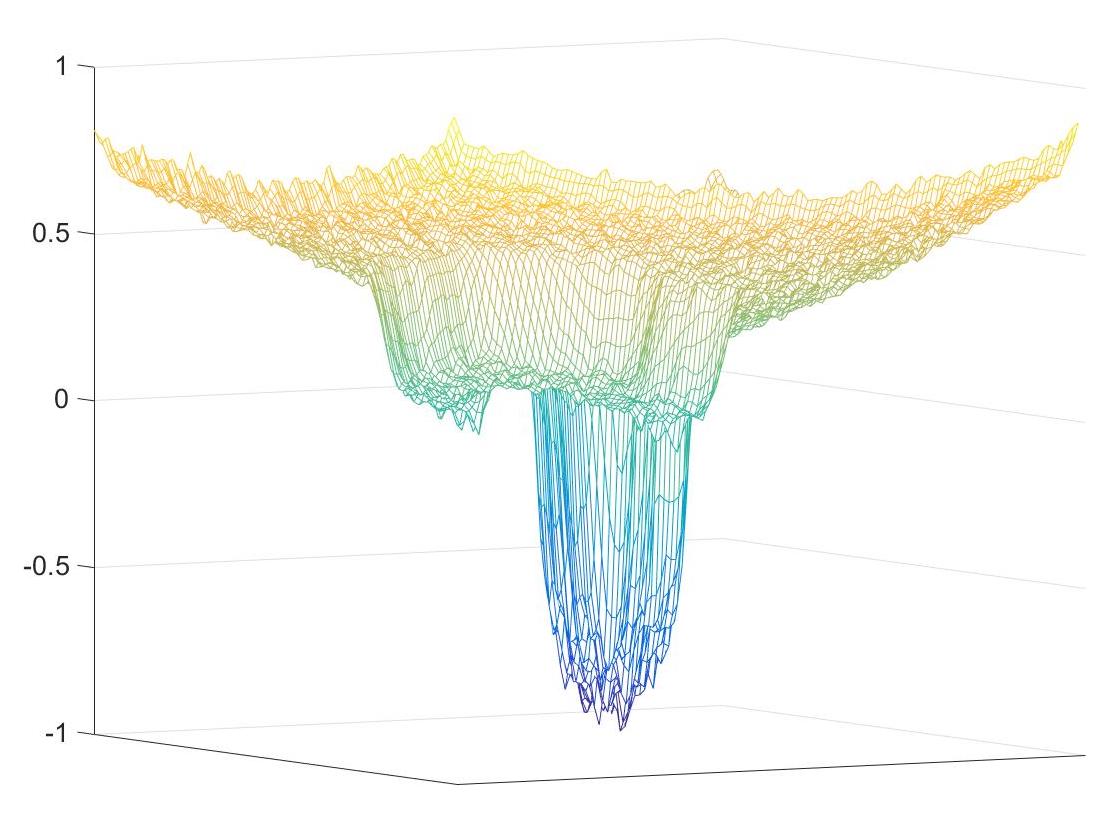}
    \end{subfigure}
    
    \begin{subfigure}{.29\linewidth}
\includegraphics[width=\linewidth]{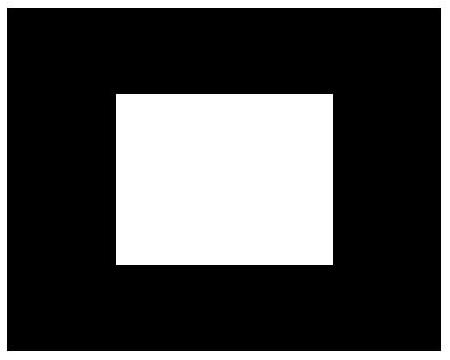}
   \end{subfigure}
    \begin{subfigure}{.29\linewidth}
\includegraphics[width=\linewidth]{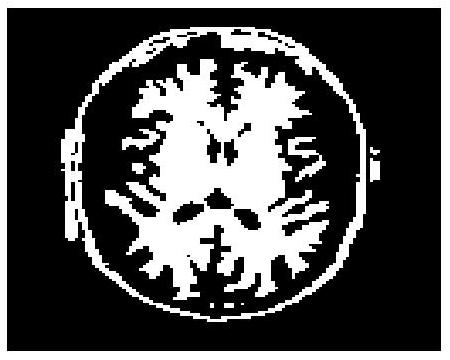}
    \end{subfigure}
    \begin{subfigure}{.29\linewidth}
\includegraphics[width=\linewidth]{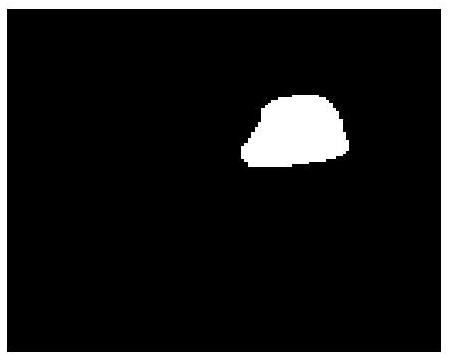}
    \end{subfigure}
\end{center}
\vspace{-20pt}
\caption{Test Problems: From left to right are Examples 1-3. From top to bottom are the image ($z$), fitting term ($f$), and thresholded segmentation result ($GT$) using the original method of \cite{Bresson:07}.} \label{fig:Examples}
\vspace{-10pt}
\end{figure}

\begin{table}[!ht]														
\begin{center}														
\begin{tabular}{| c || c | c || c | c || c | c | }
\hline	
\ \ \ \multirow{2}{*}{$q$}\ \ \ &	\multicolumn{2}{c||}{Example 1} &	\multicolumn{2}{c||}{Example 2} &	\multicolumn{2}{c|}{Example 3} \\
\cline{2-7}																	
    &	 $E_{1}$	&	 $E_{2}$	&	 $E_{1}$  	&	$E_{2}$ &	 $E_{1}$  	&	$E_{2}$		\\
\hline	 \hline																			
  0  & \ \ \ \ \ \ 0.876 \ \ \ \ \ \  &	$5.64\times 10^{2}$ & \ \ \ \ \ \ 0.943 \ \ \ \ \ \ & $2.49\times 10^{2}$ & \ \ \ \ \ \ 0.920 \ \ \ \ \ \ & $5.59\times 10^{1}$ \\	
\hline	
0.1 & 0.961 &	$1.59\times 10^{2}$ & 0.946 & $1.47\times 10^{2}$ & 0.962 & $1.19\times 10^{1}$ \\	
\hline																			
0.2 & 0.987 &	$5.16\times 10^{1}$ & 0.964 & $8.90\times 10^{1}$ & 0.969 & $1.05\times 10^{1}$ \\	
\hline																		
0.3 & 0.995 &	$1.94\times 10^{1}$ & 0.982 & $4.94\times 10^{1}$ & 0.974 & $8.02\times 10^{0}$ \\																			
\hline
0.4 & 0.999 &	$6.76\times 10^{0}$ & 0.987 & $3.47\times 10^{1}$ & 0.989 & $3.45\times 10^{0}$ \\																			
\hline		
0.5 & 1.000 &	$3.24\times 10^{0}$ & 0.988 & $2.98\times 10^{1}$ & 0.995 & $7.55\times 10^{-1}$ \\																			
\hline		
0.6 & 1.000 &	$2.65\times 10^{0}$ & 0.989 & $2.90\times 10^{1}$ & 0.999 & $6.03\times 10^{-1}$ \\																			
\hline		
0.7 & 1.000 &	$1.24\times 10^{0}$ & 0.989 & $2.90\times 10^{1}$ & 1.000 & $4.78\times 10^{-1}$ \\																			
\hline		
0.8 & 1.000 &	$5.34\times 10^{-1}$ & 0.989 & $2.90\times 10^{1}$ & 0.999 & $4.24\times 10^{-1}$ \\																			
\hline		
0.9 & 1.000 &	$2.49\times 10^{-1}$ & 0.989 & $2.73\times 10^{1}$ & 1.000 & $3.30\times 10^{-1}$ \\																			
\hline		
  1  & 1.000 & $2.10\times 10^{-5}$ & 1.000 & $3.28\times 10^{-2}$ & 1.000 & $2.10\times 10^{-1}$ \\																			
\hline																					
\end{tabular}	
\end{center}	
\vspace{-20pt}
\caption{Results: For Examples 1-3 we vary $q\in[0,1]$ and provide $E_{1}$ and $E_{2}$.} \label{tab:Main}
\vspace{-10pt}														
\end{table}	

In this section we introduce some results for the test problems shown in Figure \ref{fig:Examples}, using \eqref{eqn:CV} for $f$ in Examples 1 and 2 and \eqref{eqn:DS} for $f$ in Example 3. The focus of these results is to determine the dependence on $q$, i.e. to what extent can we restrict the domain for problems of this type? As a comparison, we use a result from \cite{Bresson:07} (for a particular choice of $\lambda$ in each case). Specifically, we iterate until the following stopping criterion is met at the $\ell^{th}$ iteration:
\begin{equation*}
\max\left\lbrace\parallel u^{(\ell)}-u^{(\ell-1)}\parallel,\parallel v^{(\ell)}-v^{(\ell-1)}\parallel\right\rbrace\leq\delta.
\end{equation*}
For $\delta=10^{-10}$ we set $u^{GT}(x)=u^{(\ell)}$ in the original dual formulation. We also use the thresholded result:
\begin{equation*}
GT(x)=\left\{
\begin{array}{ll}
1, &\text{for}\ x\in u^{GT}(x)>\epsilon \\
0, &\text{for}\ x\in u^{GT}(x)\leq\epsilon.
\end{array}
\right.
\end{equation*}
Following convention we set $\epsilon=0.5$. We refer to the solution for the proposed method (with $\delta=10^{-2}$) as $u^{*}$ and its corresponding thresholded result as $\Omega_{1}^{*}$. This allows us to define the two error measurements that we use in discussing the results when varying the parameter in the proposed method, $q$. The first is the Tanimoto Coefficient between the thresholded results, and the second is the $L^{2}$ difference between the proposed solution and the original solution:
\begin{equation*}
E_1=\frac{N(GT\cap\Omega_{1}^{*})}{N(GT\cup\Omega_{1}^{*})},\ \ \ \ \ E_2=\int_{\Omega}(u^{*}-u^{GT})^{2}\ dx.
\end{equation*}
Here $N(\cdot)$ refers to the number of nodes in the enclosed region and $E_{1}\in[0,1]$, with $E_{1}=1$ indicating a perfect result. With the second error measurement, clearly we expect $E_{2}$ to approach 0 as $q$ increases. We don't necessarily expect $E_{2}$ to be 0 for $q\leq1$ as $u^{GT}$ is not binary precisely and in the tests we use $\delta=10^{-2}$, so there is likely to be a minor difference. Whilst $E_{2}$ is a useful measure to demonstrate the correspondence between the value of $q$ and the original method, we are primarily interested in $E_{1}$ as the crucial indicator of a successful segmentation result. 

In Table \ref{tab:Main} we present the main results for $q\in[0,1]$. We include $q=0$ (i.e. a completely thresholded result) and $q=1$ (i.e. the original method) to demonstrate the full effect of the choice of $q$. Both error measurements are included (to 3 s.f.) and we can see that increasing $q$ to 0.4 is enough to produce a very good result ($E_{1}>0.98)$ in all examples. For Example 2 we can see that for $q<1$, $E_{1}$ does not reach 1 meaning that restricting the domain of the dual formulation always changes the segmentation result for this fitting term. However, there is a very close correspondence between the results even for small values of $q$, which is encouraging. As expected $E_{2}$ tends to decrease as $q$ increases, and the solution in the restricted domain is reasonably close to the original solution. To put these results in context the size of all images tested here are $128\times128$. These results demonstrate that using a dual formulation in a restricted domain is a viable approach for problems of this type.

\begin{wrapfigure}{r}{0.6\textwidth}
\vspace{-20pt}
\begin{center}
\includegraphics[width=0.6\textwidth, height=0.42\textwidth]{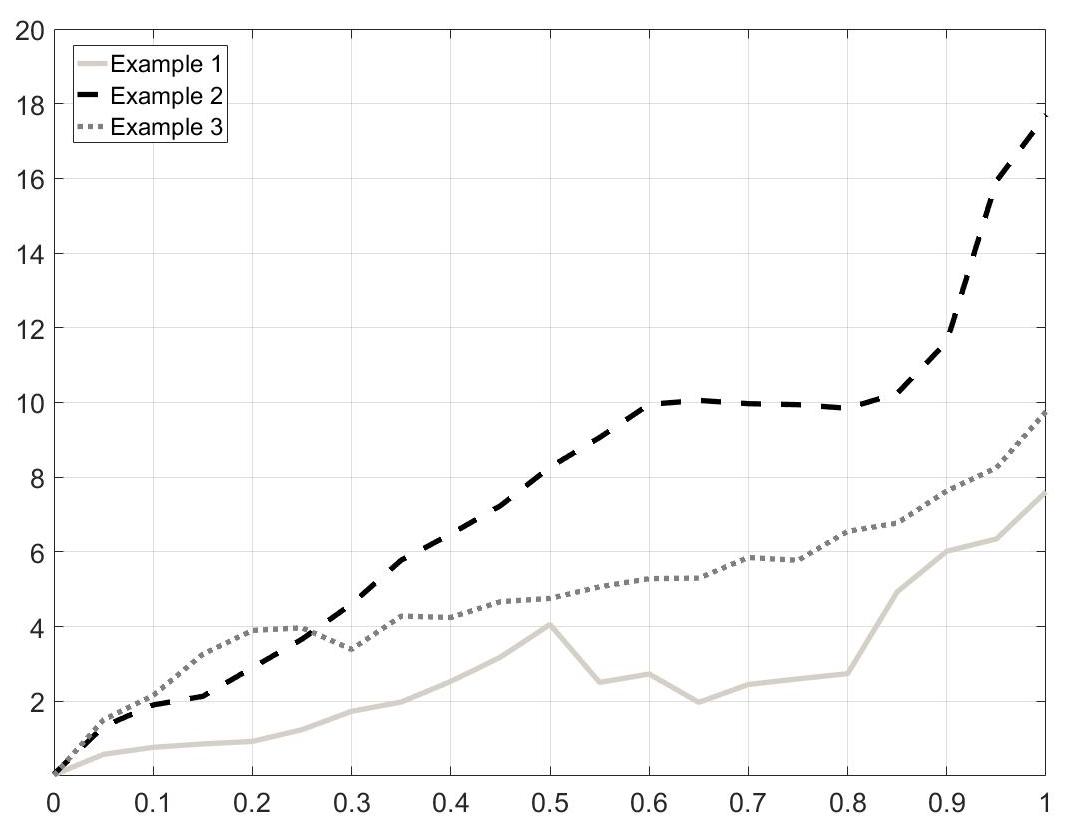}
\end{center}
\vspace{-20pt}
\caption{Computation Time, $t(q)$.} \label{fig:Time}
\vspace{-0pt}
\end{wrapfigure}

In Figure \ref{fig:Time} we include the computation time (in seconds, to 1 d.p.) for different choices of $q\in[0,1]$. Clearly the expectation is that as $q$ increases $t$ should also increase, but this does not quite hold absolutely. This is likely to be down to features in these fitting terms that mean minor increases to $q$, perhaps counterintuitively, slightly simplify the problem. For the original method of \cite{Bresson:07} the computation time was $t=6.0s$ for Example 1, $t=14.0s$ for Example 2, and $t=8.2s$ for Example 3 (which we will refer to as $t_{1},t_{2}$, and $t_{3}$ respectively). For these results, and in the following, we set $\delta=10^{-2}$ as a stopping criterion, and use them as a benchmark in each case. From Figure \ref{fig:Time} it can be seen that for Examples 1-3 $t$ is only greater than $t_{1},t_{2}$, or $t_{3}$ for $q>0.9$. For smaller values of $q$ it is possible to make significant gains in terms of computation time. For Example 1, when $E_{1}>0.98$ the average time is $t=2.8s$ for $q\leq0.9$. Similarly for Examples 2 and 3, the average times are $t=8.7s$ and $t=5.5s$ respectively. This corresponds to a time saving of $53\%$, $38\%$, and $33\%$ for Examples 1-3, for cases with an accurate segmentation. Extending these tests to a wider choice of fitting terms, and investigating the effect of changing $\lambda$ and $\delta$, would help further determine the effectiveness of the proposed approach.

\begin{figure}[!ht]
\begin{center} 
    \begin{subfigure}{.28\linewidth}
\includegraphics[width=\linewidth]{Hf1.jpg}
   \end{subfigure}
    \begin{subfigure}{.28\linewidth}
\includegraphics[width=\linewidth]{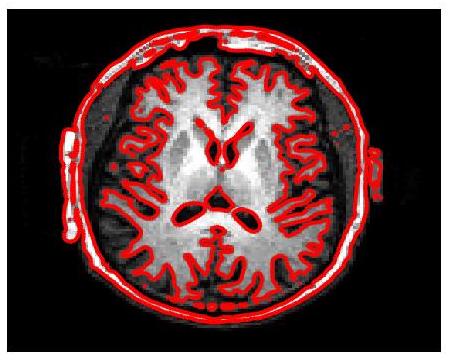}
    \end{subfigure}
    \begin{subfigure}{.28\linewidth}
\includegraphics[width=\linewidth]{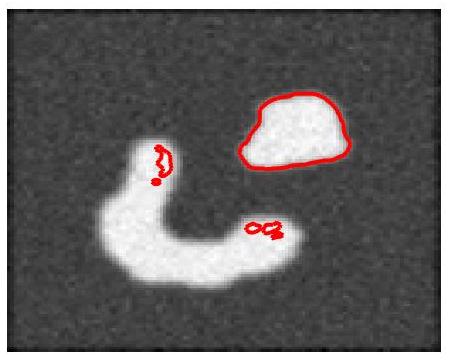}
    \end{subfigure}

    \begin{subfigure}{.28\linewidth}
\includegraphics[width=\linewidth]{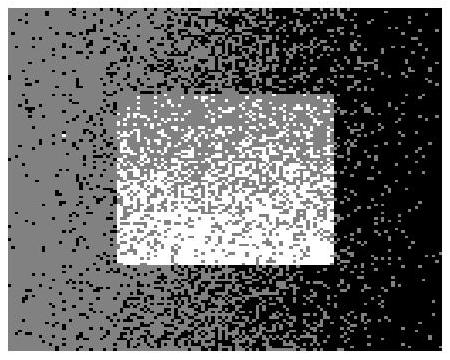}
   \end{subfigure}
    \begin{subfigure}{.28\linewidth}
\includegraphics[width=\linewidth]{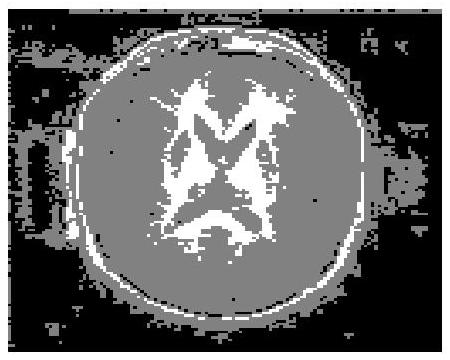}
    \end{subfigure}
    \begin{subfigure}{.28\linewidth}
\includegraphics[width=\linewidth]{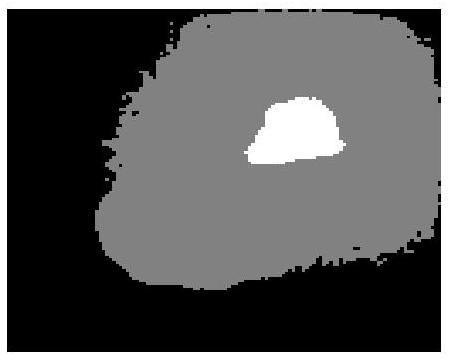}
    \end{subfigure}
    
   \begin{subfigure}{.28\linewidth}
\includegraphics[width=\linewidth]{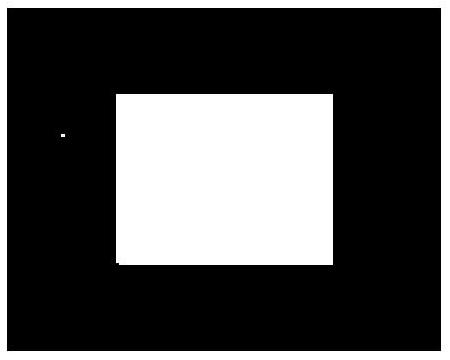}
   \end{subfigure}
    \begin{subfigure}{.28\linewidth}
\includegraphics[width=\linewidth]{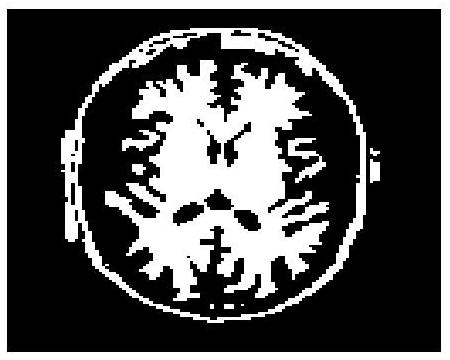}
    \end{subfigure}
    \begin{subfigure}{.28\linewidth}
\includegraphics[width=\linewidth]{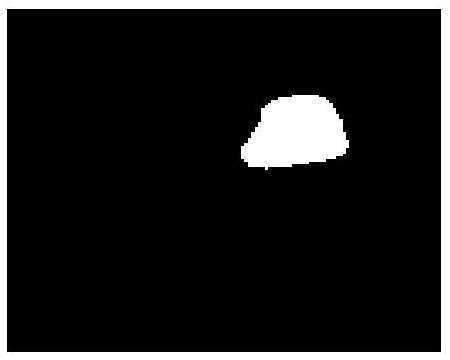}
    \end{subfigure}
\end{center}
\vspace{-20pt}
\caption{Results: From left to right are Examples 1-3, respectively. From top to bottom are the image with the zero level-set of $f$ in red, the regions $Fg$ (white), $Bg$ (black), $RD$ (grey), and $\Omega_{1}^{*}$ for $q=0.5$.} \label{fig:Results}
\vspace{-10pt}
\end{figure}

In Figure \ref{fig:Results} we present some example results for $q=0.5$. For Example 1, $E_{1}=1.000$ and $t=4.1s$. For Example 2, $E_{1}=0.988$ and $t=8.3s$. For Example 3, $E_{1}=0.995$ and $t=4.8s$. Compared to the original method, $t_{1}=6.0s$, $t_{2}=14.0s$, and $t_{3}=8.2s$. In each case a significant improvement can be made in terms of computation time with minimal compromise on the quality of the result as measured by $E_{1}$.

\section{Conclusion}

The results presented support the idea that the domain can be restricted for problems of this type, without compromising the quality of the result. This allows for significant gains in terms of computation time. Additional testing to verify these findings would be beneficial, particularly for a wider variety of fitting terms. An example would be where $f$ contains high levels of noise. Further considerations might be necessary to restrict the domain in a robust way, and developing the required understanding would help develop this approach further. For $q=1$ the computation times for Examples 1-3 are $t=7.6s, 17.7s$, and  $9.8s$ respectively, which corresponds to approximately a $25\%$ increase for the proposed method when no restriction of the domain is considered. If the efficiency of the domain restriction could be improved this would allow for higher values of $q$ to be selected for a reduced cost, which could be particularly beneficial in cases of high noise in the fitting term.

The results presented are for images of size $128\times128$ in order to explore the viability of restricting the domain for this problem. Improvements in $t$ for larger images, or 3D problems, could be particularly valuable. Extending this approach to these cases is of interest, and would help support the proposed approach further. We are also considering the extension of this approach beyond the dual formulation of \cite{Bresson:07}, such as split-Bregman \cite{Goldstein:10} and additive operator splitting \cite{CDSS}. Following the framework introduced here, assumptions about the solution of \eqref{eqn:CR} can be adapted to other methods in a similar way. However, the initial results presented here are encouraging.

\section*{Acknowledgements}

The author would like to acknowledge the support of the EPSRC grant EP/N014499/1.


\bibliographystyle{apalike}

\bibliography{imvip2017}

\end{document}